\documentclass{IEEEtran}
\usepackage{isolatin1}
\usepackage{times}
\usepackage{pdfsync}
\usepackage{epsfig}
\usepackage{epstopdf}
\usepackage{amsmath}
\usepackage{amssymb}
\usepackage{mathrsfs}
\usepackage{color}
\usepackage{fix2col}

 \newtheorem{theorem}{Theorem}[section]
 \newtheorem{lemma}[theorem]{Lemma}

 \newtheorem{remark}[theorem]{Remark}

 \newtheorem{problem}[theorem]{Problem}

\title{Transmission Power Scheduling for Energy Harvesting Sensor in Remote State Estimation$^\ast$}
\author{Yuzhe Li, Daniel E.\ Quevedo, Vincent Lau, Subhrakanti Dey, and Ling
  Shi%
\thanks{Yuzhe Li, Vincent Lau and Ling Shi are with the Department of
      Electronic and Computer Engineering, Hong Kong University of Science and
      Technology, Clear Water Bay, Kowloon, Hong Kong (emails: \{yliah, eeknlau,
      eesling\}@ust.hk). Daniel Quevedo is  with the School of
    Electrical Engineering \& 
      Computer Science, The University of Newcastle, NSW
      2308, Australia; e-mail:
      dquevedo@ieee.org. Subhrakanti Dey is with 
the Department of Engineering Sciences, Uppsala University, Sweden
(email:subhrakanti.dey@angstrom.uu.se).}\thanks{This research was supported
      by an HK RGC GRF grant 618612 and under Australian Research Council's 
  Discovery Projects funding scheme (project number DP0988601).}
\thanks{$^\ast$ Extended version of article to be published in Proc.\ 19th IFAC World Congress 2014.}}

\begin{document}
\maketitle

\begin{abstract}
We study remote estimation in a wireless sensor network. Instead of using a conventional battery-powered sensor, a sensor equipped with an energy harvester which can obtain energy from the external environment is utilized. We formulate this problem into an infinite time-horizon Markov decision process and provide the optimal sensor transmission power control strategy. In addition, a sub-optimal strategy which is easier to implement and requires less computation is presented. A numerical example is provided to illustrate the implementation of the sub-optimal policy and evaluation of its estimation performance.
\end{abstract}

%
%
%
%


\section{Introduction}

Wireless sensors network (WSN) has been a hot research topic in recent years. Both theoretical results and practical applications are growing rapidly.  Compared with traditional wired sensors, wireless sensors provide many advantages such as low cost, easy installation, and self-power. In a WSN, sensors are typically equipped with batteries and expected to work for a long time (\cite{yick2008wireless}). Thus, the energy constraint is an inevitable issue.  In some applications, the amounts of sensors can be quite large (e.g., environment monitoring) or sensors may be located in dangerous environments (\cite{ho2012optimal}) (e.g., chemical industry), making the replacement of batteries difficult or even impossible.

To deal with energy aspects of WSN, one possible way is to develop more efficient sensor energy power control methods to make the best use of the batteries (\cite{aziz2013survey,pantazis2007survey,quevedo2010energy}). Those existing results demonstrate significant improvement of the lifetime of the sensor and system performance under energy constraints. The problem is, however, still not completely solved as the the battery will eventually run out. At the same time, the optimization of lifetime of the sensor under limited energy will always lead to other sacrifices such as estimation quality or system stability (\cite{sudevalayam2011energy}).

To overcome this limitation, an alternative way is to replace the conventional
battery-powered sensor with sensors equipped with an energy harvester. The
technology of energy harvesting refers to obtaining energy from the external
environment or other types of energy sources (e.g., body heat, solar energy,
piezoelectric energy, wind energy) and converting them into electrical energy
which can be stored and used by the sensor  (\cite{ho2012optimal}). For sensors
using this technology, the energy (but not the energy-rate) is typically
``unlimited'' compared to battery-powered sensor as the harvester can generate
power all the time during the whole time-horizon. But unlike the battery-powered
sensor which has relatively explicit energy amount for future use, the sensor
with energy harvester will be subject to an unpredictable future energy level as
they are affected by the external environment. Due to the randomness of the
amounts of harvested energy in the following time steps, new challenges arise in
the design and analysis of the communication strategy of the sensor. Power
control and battery management requires trading off current transmission success
probabilities for expected future ones. 

The work \cite{ho2012optimal} studied the problem of energy allocation for wireless communication. The authors aimed to maximize the throughput under time-variant channel conditions and harvested energy sources, which is solved by dynamic programming and convex optimization techniques. In \cite{nayyar2012optimal}, the authors investigated a remote estimation problem for an energy harvesting sensor and a remote estimator. The communication strategy for the sensor and the estimation strategy for the remote estimator are jointly optimized in terms of the expected sum of communication and distortion costs, again using a dynamic programming approach.

In our preliminary work, \cite{li2013optimal}, an optimal periodic sensor power schedule is derived. The proposed method is, however, only suitable for solving the problem of battery-powered sensor subject to an average energy constraint. For energy harvesting sensor, a new approach is needed to handle the randomness of the energy constraints. Driven by this motivation, in the present work, we consider remote estimation with a wireless sensor having an energy harvesting capability. The most related result of our present work is \cite{norican13optimal}, which studied optimal transmission energy allocation scheme for error covariance minimization in Kalman filtering with random packet losses when the sensors have energy harvesting capabilities, and they provided some structural results on the optimal solution for both finite and infinite time-horizon. Different from their work, we specify the different distributions of different environment conditions for the energy harvesting model. Furthermore, we use a smart sensor to pre-processes the measurement data which can improve the estimation quality \cite{hovareshti2007sensor}. The main challenges and contributions of this work are summarized as follows:
\begin{enumerate}
\item \textbf{Randomness of harvested energy:} In previous works, e.g., \cite{li2013optimal,quevedo2010energy}, the constraints of the transmission power are deterministic. For energy harvesting sensors, on the other hand, the information of the energy constraints is not exactly available for the sensor before the harvesting due to the randomness of the energy resources. To handle this new challenge, we  develop a new approach.

\item \textbf{Infinite time-horizon MDP:} We consider an infinite time-horizon problem, which is a better approximation for long-run applications and more difficult. In order to overcome the randomness of the energy resources, we prove that an associated power control design problem can be formulated into a standard MDP framework with infinite time-horizon and give the optimal solution.

\item  \textbf{Sub-optimal solution:} As the MDP method cannot in general provide an explicit form of the optimal solution and the computational complexity is formidable for general higher-order systems, we propose a sub-optimal solution which is in threshold form and is easy to implement for different system parameter settings.
\end{enumerate}

The remainder of this manuscript is organized as follows. Section 2 presents the system setup. Section 3 formulates the problem into a standard MDP framework and provides the optimal solution. Section 4 introduces a sub-optimal solution and compares it with the optimal one. Numerical example and simulations are included in Section 5. Section 6 draws conclusions.

\textit{Notations}: $\mathbb{Z}$ denotes the set of integers and $\mathbb{N}$ the positive integers. $\mathbb{R}$ is the set of real numbers. $\mathbb{R}^{n}$ is the $n$-dimensional Euclidean space. $\mathbb{S}_{+}^{n}$ (and $\mathbb{S}_{++}^{n}$) is the set of $n$ by $n$ positive semi-definite matrices (and positive definite matrices). When $X \in\mathbb{S}_{+}^{n}$ (and $\mathbb{S}_{++}^{n}$), we write $X \geqslant 0$ (and $X > 0$). $X\geqslant Y$ if $X - Y \in \mathbb{S}_{+}^{n}$. $\mathrm{Tr}(\cdot)$ is the trace of a matrix. The superscript $'$ stands for transposition. For functions $f, f_1, f_2$ with appropriate domains, $f_1\circ f_2(x)$ stands for the function composition $f_1\big(f_2(x)\big)$, and $f^{n}(x) \triangleq f\big(f^{n-1}(x)\big)$, where $n\in\mathbb{N}$ and with $f^{0}(x) \triangleq x$. $\delta_{ij}$ is Dirac delta function, i.e., $\delta_{ij}$ equals to $1$ when $i=j$ and $0$ otherwise. The notation $\mathbb{P}[\cdot]$ refers to probability and $\mathbb{E}[\cdot]$ to expectation.

\section{State Estimation With An Energy Harvester } \label{sec:Pre}

We consider the problem of remote estimating the state of the following linear time-invariant (LTI) system:
\begin{align}
  x_{k+1} & =  Ax_k + w_k, \label{eqn:process-dynamics} \\
  y_k & = Cx_k + v_k,  \label{eqn:measurement-equation}
\end{align}
where $k\in \mathbb{N}$, $x_k \in \mathbb{R}^{n_x}$ is the system state vector at time $k$, $y_k\in \mathbb{R}^{n_y}$ is the measurement taken by the sensor, $w_{k} \in\mathbb{R}^{n_x} $ and $v_k \in \mathbb{R}^{n_y}$ are zero-mean i.i.d. Gaussian noises with $\mathbb{E}[w_{k}w_{j}^\prime] =\delta_{kj}Q$ ($Q\geqslant 0$), $\mathbb{E}[v_{k}(v_{j})^\prime] = \delta_{kj}R$ ($R > 0$), $\mathbb{E}[w_{k}(v_{j})^\prime] = 0 \; \forall j,k\in\mathbb{N}$. The initial state $x_0$ is a zero-mean Gaussian random vector with covariance $\Pi_0\geqslant 0$ and is uncorrelated with $w_k$ and $v_k$. The pair $(A, C)$ is assumed to be observable and $(A, Q^{1/2})$ is controllable.

\begin{figure}[htp]
  \centering
  \includegraphics[width=9cm]{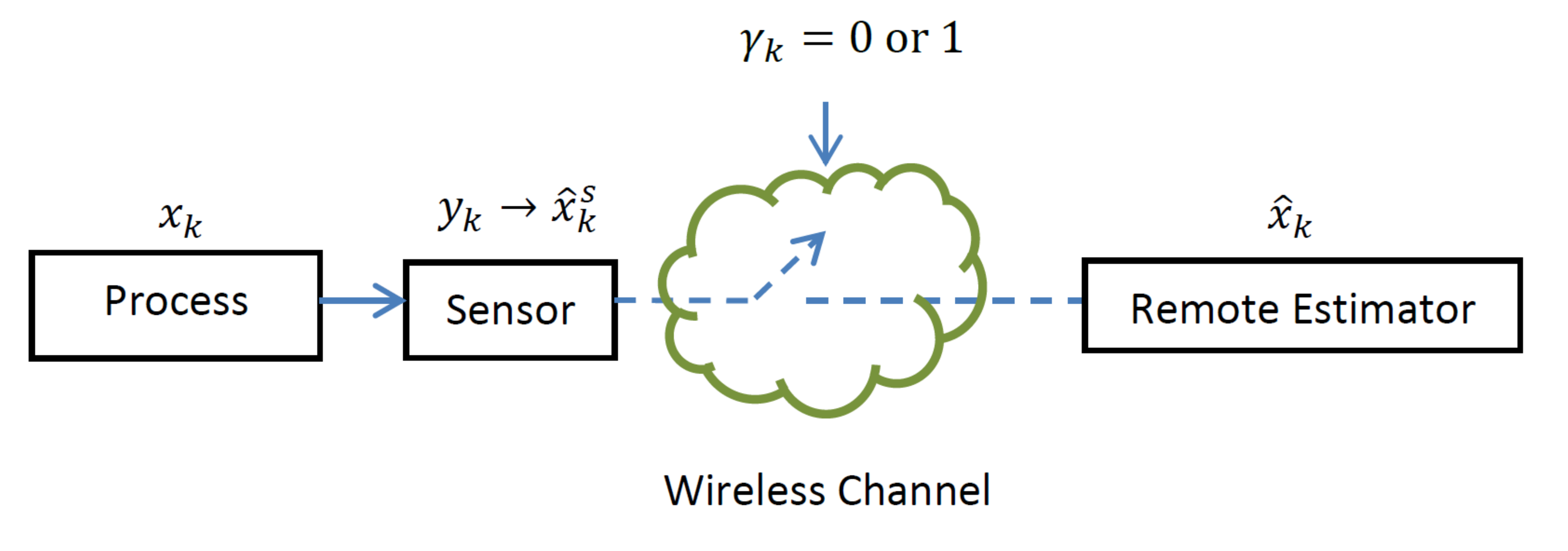}
  \caption{System Architecture} \label{fig:system}
\end{figure}

\subsection{Sensor Local State Estimate}
We assume the sensor in this work is embedded with an on-board processor (\cite{hovareshti2007sensor}), the so called ``smart sensor'' (\cite{shi2012scheduling,wu2013can,li2013jamming}). At each time $k$, the sensor first locally runs a regular Kalman filter to produce the minimum mean-square error (MMSE) estimate of the state $x_k$ based on all the measurements it collects up to time $k$. It then transmits the local estimate to a remote estimator.

Denote $\hat{x}_k^s$ and $P_k^s$ as the sensor's local MMSE state estimate and the corresponding estimation error covariance, respectively, i.e.:
\begin{align} \label{eqn:error-covariance-all}
  \hat{x}_k^s&=\mathbb{E}[x_k|y_1,y_2,...,y_k],\\
  P_k^s&=\mathbb{E}[(x_k-\hat{x}_k^s)(x_k-\hat{x}_k^s)'|y_1,y_2,...,y_k],
\end{align}
which can be calculated recursively using standard Kalman filter update equations (\cite{anderson1981detectability}):
\begin{align}
  \hat{x}_{k|k-1}^s  = & A\hat{x}_{k-1}^s, \label{eqn:state-time-update} \\
  P_{k|k-1}^{s}  = & AP_{k-1}^sA' + Q, \label{eqn:error-covariance-time-update} \\
  K_{k}^s  = & P_{k|k-1}^{s}C'[CP_{k|k-1}^{s}C' + R]^{-1}, \label{eqn:kalman-gain} \\
  \hat{x}_{k}^s  = & A\hat{x}_{k-1}^s + K_{k}^s(y_{k} - CA\hat{x}_{k-1}^s), \label{eqn:state-measurement-update} \\
  P_{k}^s  = &(I -K_{k}^sC)P_{k|k-1}^{s}, \label{eqn:error-covariance-measurement-update}
\end{align}
where the recursion starts from $\hat{x}_{0}^s = 0$ and $P_0^s = \Pi_0\geqslant0$.

The following Lyapunov and Riccati operators $h,\tilde{g}: \mathbb{S}_{+}^{n}\rightarrow \mathbb{S}_{+}^{n}$ are introduced to facilitate our subsequent discussion:
\begin{align}
h(X) &\triangleq  AXA' + Q,\label{eqn:definition_h}\\
\tilde{g}(X) &\triangleq  X - XC'[CXC' + R]^{-1}CX.\label{eqn:definition_tilde_g}
\end{align}

Since the estimation error covariance $P_k^s$ in (\ref{eqn:error-covariance-measurement-update}) converges to a steady-state value exponentially fast (See \cite{anderson1981detectability}), without loss of generality, we assume that the Kalman filter at the sensor side has already entered the steady state, i.e., :
\begin{equation}\label{eqn:steady-state-assumption}
  P_k^s = \overline{P},~k \geqslant 1,
\end{equation}
where $\overline {P}$ is the steady-state error covariance, which is the unique positive semi-definite solution of $\tilde{g}\circ h(X) = X$.

$\overline {P}$ has the following property (see \cite{shi2011time}).

\begin{lemma} \label{lemma:steady-state-error-property}
For $0\leqslant t_1\leqslant t_2$, the following inequality holds (\cite{shi2011time}):
\begin{equation}
    h^{t_1}(\overline{P}) \leqslant h^{t_2}(\overline{P}).
\end{equation}
In addition, if $t_1 < t_2$, then
\begin{equation} \label{eqn:steady-state-error-property}
    \mathrm{Tr}\left(h^{t_1}(\overline{P})\right) < \mathrm{Tr}\left(h^{t_2}(\overline{P})\right).
\end{equation}
\end{lemma}

\subsection{Wireless Communication Model}

Similar to~\cite{li2013online}, c.f.,\cite{quevedo2010energy},  the local state estimate of the sensor $\hat{x}^s_k$ is transmitted to the remote estimator over an Additive White Gaussian Noise (AWGN) channel using Quadrature Amplitude Modulation (QAM). Denote $\omega_k$ as the transmission power for sending the QAM symbol at time $k$, which will be designed in the following sections. Based on the analysis in \cite{li2013online},  the approximate relationship between the symbol error rate (SER) and $\omega_k$ is given by
\begin{equation} \label{eqn:SER1}
  \text{SER}\approx\exp\Big(-\beta \frac{\omega_k}{N_0W}\Big).
\end{equation}

The communication channel is assumed to be time-invariant, i.e., $\beta$, $N_0$, $W$, are constants during the whole time horizon \footnote{For time-variant channels, one can also formulate the problem in a similar way. This is left for future work}. In practice, the remote estimator can detect symbol errors via cyclic redundancy check (CRC). Thus taking into account of the SER in the transmission of QAM symbols, a binary random process $\{\gamma_k\}, k\in\mathbb{N}$ can be used to characterize the equivalent communication channel for $\hat{x}_k^s$ between the sensor and the remote estimator, where
\begin{equation} \label{eq:43}
  \gamma_k =
  \begin{cases}
  1, &\text{if $\hat{x}^s_k$ arrives error-free at time $k$,}\\
  0, & \text{otherwise (regarded as dropout).}
  \end{cases}
\end{equation}

From (\ref{eqn:SER1}), we have
\begin{equation} \label{eqn:drop model-1}
\mathbb P[\gamma_k=0]=(1-\lambda)^{\omega_k},
\end{equation}
where $\lambda$ is given by:
\begin{equation}\label{def:lambda}
  \lambda\triangleq 1-\exp( -\frac{\beta}{N_0W})\in(0,1).
\end{equation}

%

\subsection{Energy Harvester}
Now we present a simple model for  the energy harvesting sensor.  Assume that there are two states of the external environment: $G$ denotes the good condition (e.g., windy, sunny, etc.) and $B$ denotes the bad condition which may alternate at every time step. At time $k$, the environment condition state is denoted as $e_k$ and the transition of the two condition states between two time steps follows a Markov chain model:
\begin{figure}[t]
  \centering
  \includegraphics[width=8.5cm]{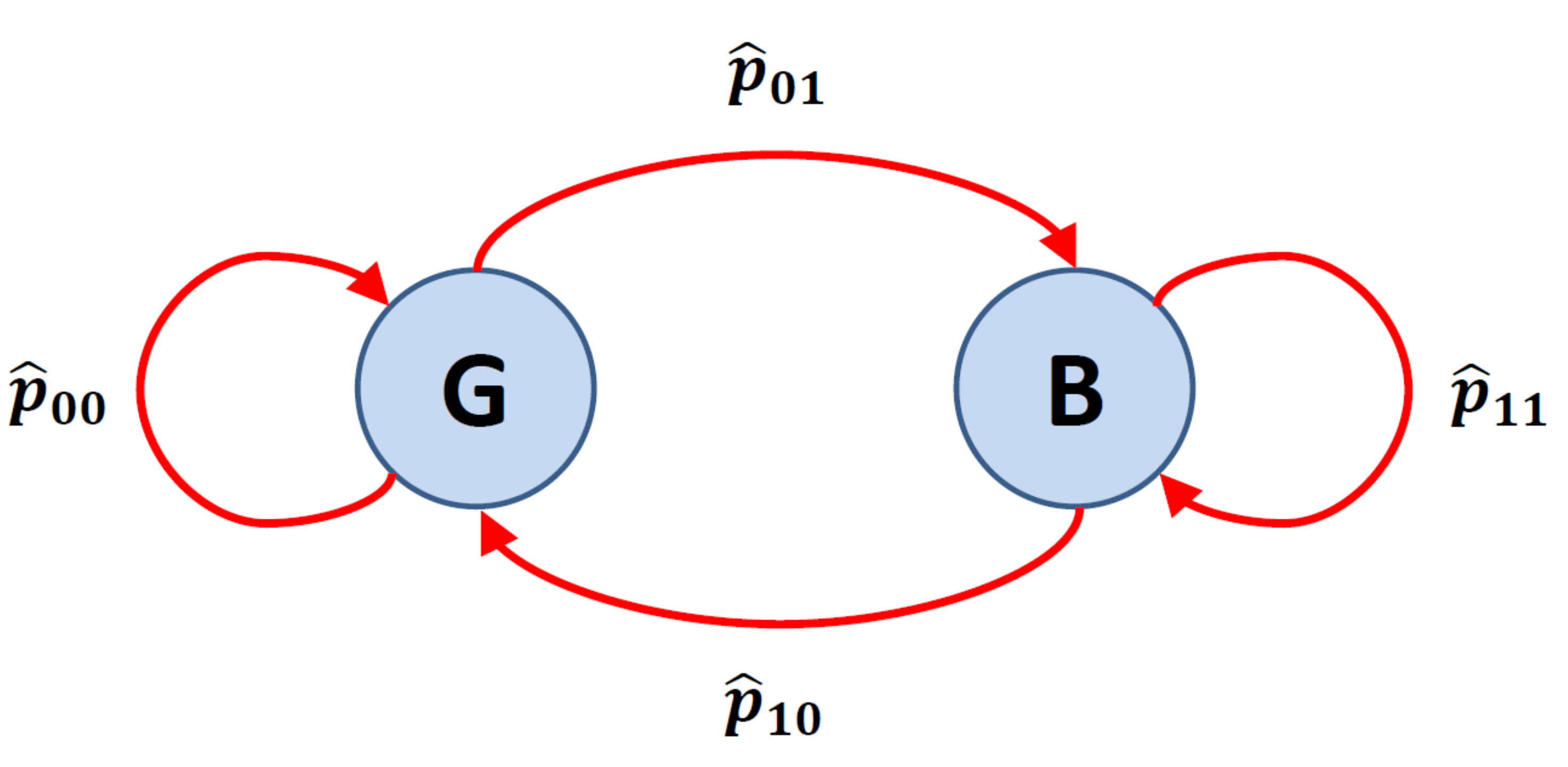}\\
  \caption{Markov Chain Model of Environment Condition}
\end{figure}

The transition can be expressed as
\begin{align}\label{eqn:environ_trans}
\mathbb P(e_{k+1}=G|e_k=G)&=\hat p_{00} ,\\
\mathbb P(e_{k+1}=B|e_k=G)&=\hat p_{01},\\
\mathbb P(e_{k+1}=G|e_k=B)&=\hat p_{10} ,\\
\mathbb P(e_{k+1}=B|e_k=B)&=\hat p_{11}.
\end{align}

Denote the remaining energy level in the sensor's battery at the beginning of time step $k$ as $b_k$. The maximum battery level (battery capacity) is denoted as $b_{\text{max}}$. At each time step, we assume the amount of harvested energy, denoted as $r_k$, is a discrete random variable which can only take values in $\mathbb N^+$, i.e., $r_k\in\{0,1,2,..., b_{\text{max}}\}$ (Note that for the situation that $r_k>b_{\text{max}}$, we can regard it as $r_k=b_{\text{max}}$ and add up all the corresponding probabilities as $\mathbb P[r_k=b_{\text{max}}]$). Under different environment conditions, $r_k$ follows different distributions:
\begin{equation} \label{eqn:harvested_prob1}
\mathbb P[r_k=i| e_k=G]=\pi_i^0,
\end{equation}
and
\begin{equation} \label{eqn:harvested_prob2}
\mathbb P[r_k=i| e_k=B]=\pi_i^1,
\end{equation}
where $i\in\{0,1,2,..., b_{\text{max}}\}.$

Note that after harvesting the energy $r_k$, the battery level now is  $\min\{b_k+r_k, b_{\rm{max}}\}$. Then the sensor needs to decide the transmission power $\omega_k$ used at time $k$ to send the local state estimates to the remote estimator based on the current battery level. After this procedure, the process moves to next time step $k+1$ and the battery level at the beginning of $k+1$ is
\begin{equation}\label{eqn:battery_recur}
b_{k+1}=\min\{b_k+r_k, b_{\rm{max}}\}-\omega_k.
\end{equation}

As mentioned before, different power levels lead to different dropout rates, and thereby affect the estimation performance. Whilst keeping the battery partly charged serves to ``prepare for the future'', one should also avoid wasting energy harvesting opportunities due to the battery being full.  This motivates the issue of energy management to be studied in Section III.

\subsection{Remote State Estimation}

Denote $\hat{x}_k$ and $P_k$ as the remote estimator's own MMSE state estimate and the corresponding error covariance based on all the  sensor data packets received up to time step $k$. The works \cite{li2013optimal} and ~\cite{shi2010kalman} show that they can be calculated via the following procedure: once the sensor's local estimate arrives, the estimator synchronizes $\hat{x}_k$ with that of the sensor, i.e., $\hat{x}_k^s$; otherwise, the remote estimator just predicts $\hat{x}_k$ based on its previous estimate using the system model \eqref{eqn:process-dynamics}. From \eqref{eq:43}, the remote state estimate $\hat{x}_k$  thus obeys the recursion
\begin{equation}\label{eqn:KF-remote-estimator}
\hat{x}_k = \left\{\begin{array}{ll}\hat{x}_{k}^s, &  \mathrm{if}~\gamma_k=1,\\
A\hat{x}_{k-1},  & \mathrm{if}~\gamma_k=0.
\end{array}\right.
\end{equation}
The corresponding state estimation error covariance $P_k$ satisfies
\begin{equation}\label{eqn:KF-remote-error-covariance}
P_k = \left\{\begin{array}{ll}\overline{P}, &  \mathrm{if}~\gamma_k=1,\\
h(P_{k-1}),  & \mathrm{if}~\gamma_k=0.
\end{array}\right.
\end{equation}

\section{Optimal Transmission Power Schedule} \label{sec:MDP}

The objective of the remote estimator is to give accurate state estimates $\hat x_k$. To be more specific, we seek to minimize the trace of the average expected state estimation error covariance:
\begin{align} \label{eqn:objective-func}
  J(\theta)=\limsup_{T\to\infty}\frac{1}{T}\sum_{k=1}^{T}\mathrm{Tr}
  \{\mathbb{E}[P_k]\},
\end{align}
where $\theta=\{\omega_1, \omega_2, ...\}$ is the transmission power used at each time step. Note that here we consider an infinite time-horizon.

Due to the energy and battery constraints, we are interested in finding the optimal transmission power policy $\theta^\star$ for the sensor that solves the following constrained optimization problem:
\begin{problem} \label{problem:problem-1}
\begin{align*}
    &\min_{\theta}~~~J(\theta)\\
    &\mathrm{s.t.} ~~~0\leqslant\omega_k\leqslant \min\{b_k+r_k, b_{\rm{max}} \},~~\forall k\in \mathbb N^+,
\end{align*}
where $\theta=\{\omega_1, \omega_2, ...\}.$
$\hfill \blacksquare$
\end{problem}

We will next formulate the optimization  in Problem \ref{problem:problem-1} as an MDP problem and study the optimal policy.

As described before, the amount of the harvested energy is discrete. For convenience, we assume that the sensor can choose transmission power discretely, i.e.,
\begin{equation*}
\omega_k\in\Big\{0,1,2,...,\min\{b_k+r_k, b_{\rm{max}}\}\Big\}.
\end{equation*}

We assume that the remote estimator will send ACKs to the sensor to indicate whether it has received the data packet successfully or not (\cite{li2013online}) at time $k$, which enables the sensor to obtain $P_{k-1}$. Accordingly we define the state for the power management problem at the beginning of time step $k$ as:
\begin{equation*}
\Phi_k=\Big(\min\{b_k+r_k, b_{\rm{max}}\}, ~E_k, ~P_{k-1}\Big),
\end{equation*}
which consists of the battery level at the beginning of time step $k$:
\begin{equation*}
 \min\{b_k+r_k, b_{\rm{max}}\},
\end{equation*}
the environment condition:
\begin{equation*}
  E_k =
  \begin{cases}
  1, &\text{if}~ e_k=G,\\
  0, & \text{if}~ e_k=B,\\
  \end{cases}
\end{equation*}
and the state estimate error covariance of the last time step $P_{k-1}$. Note that here we choose $P_{k-1}$ because $P_k$ is still unknown at the beginning of time step $k$. The initial state is denoted as $\Phi_0=\phi_0$.

\begin{remark}
  $\hfill \blacksquare$
\end{remark}

From the recursion of $P_k$ in~(\ref{eqn:KF-remote-error-covariance}), it is easy to see that at any time step $k_2\geqslant k_1$, $P_{k_2}$ can be written as $P_{k_2}=h^{k_2-k_1}(\overline {P})$, where $k_1$ is the latest time when it successfully received sensor data. Since $P_k$ only takes value in the set of $\{\overline P, h(\overline P), h^2(\overline P),...\}$, the state space $\mathbb S$ for $\Phi_k$ is countably infinite:
\begin{equation*}
\mathbb S=\big\{(m,n,l)\big\},
\end{equation*}
where
\begin{align*}
  &m\in\{0,1,...,b_{\text{max}}\},\\
  &n\in\{0,1\},\\
  &l\in\{\overline P, h(\overline P),...\}.
\end{align*}

At each time step $k$, the action for the remote estimator is defined as the transmission power $\omega_k$ it chooses. Thus the available actions set $\mathbb A_k$ for time step $k$ is also finite:
\begin{equation*}
\mathbb A_k=\big\{0,1,2,..., \min\{b_k+r_k, b_{\rm{max}}\} \big\},
\end{equation*}
and therefore the action set $\mathbb A$ is
\begin{equation*}
\mathbb A\triangleq \bigcup_{k=1}^{+\infty} \mathbb A_k= \big\{0,1,2,...,  b_{\rm{max}} \big\}.
\end{equation*}

From  Section \ref{sec:Pre}, it is easy to show that the random process $\Phi:=\{\Phi_k\}$ combined with the action $\{\omega_k\}$  constitute an MDP \cite{puterman2009markov}. Define the transition probabilities $\mathbb T: \mathbb S\times \mathbb A\rightarrow \mathbb P[\mathbb S]$ as the description of each action's effect in the next state and :
\begin{equation*}
p_k( \phi_2| \phi_1, a)=\mathbb P(\Phi_{k+1}=\phi_2|\Phi_{k}=\phi_1, \omega_k=a),~~\forall k\in\mathbb N^+.
\end{equation*}

As the functions $p_k( \phi_2| \phi_1, a)$ do not depend on $k$, i.e., $\Phi$ is a time-homogeneous process, we can write $p( \phi_2| \phi_1, a)$ instead of $p_k( \phi_2| \phi_1, a)$:
\begin{equation*}
p( \phi_2| \phi_1, a)=\mathbb P(\Phi_{k+1}=\phi_2|\Phi_{k}=\phi_1,\omega_k=a),~~\forall k\in\mathbb N^+.
\end{equation*}

The closed-form expression for the one-step transition probabilities can be derived as follows.

Assume that at time $k$, the state is $\Phi_{k}=(m,n,l)$, i.e., the remaining battery level at the beginning of this time step is $m$, the environment condition index is $n$ ($n=0$ denotes good condition and $n=1$ denotes bad condition), and $P_{k-1}=l$. Though $P_{k}$ can take value from a countably infinite set, once $P_{k-1}$ is given, based on the recursion in \eqref{eqn:KF-remote-error-covariance}, there are only two possible states for $P_{k}$: $h(P_{k-1})$ and $\overline P$, with probability $(1-\lambda)^{\omega_k}$ and $1-(1-\lambda)^{\omega_k}$, respectively.  After the sensor chooses the transmission power $\omega_k$ and sends the data packet carrying $\hat x_k^s$, we can calculate the probabilities of different values $P_k$ may take.

Suppose that
\begin{equation*}
\Phi_{k+1}=\Big(m',n',l'\Big),
\end{equation*}
where $l'=P_k$ and
\begin{equation*}
m'=\min\{b_{k+1}+r_{k+1}, b_{\rm{max}}\}.
\end{equation*}

Since
\begin{align*}
b_{k+1}&=\min\{b_k+r_k, b_{\rm{max}}\}-\omega_k,\\
&=m-\omega_k,
\end{align*}
we also have
\begin{equation*}
m'=\min\{m-\omega_k+r_{k+1}, b_{\rm{max}}\}.
\end{equation*}

Clearly, when $l'\neq \overline P$ or $h(l)$,  we have
\begin{equation*}
  p\big((m', n',l')|(m,n,l), \omega_k\big)=0.
\end{equation*}

Based on the battery level recursion in \eqref{eqn:battery_recur} and environment condition transition in \eqref{eqn:environ_trans}, when $m'<b_{\rm{max}}$,  indicating that $r_{k+1}=m'- (m-\omega_k)<b_{\rm{max}}-(m-\omega_k)$, we have
\begin{equation*}
p\Big(\big(m', n',h(l)\big)|(m,n,l),\omega_k\Big)=(1-\lambda)^{\omega_k}\hat p_{nn'}\pi^{n'}_{r_{k+1}},
\end{equation*}
and
\begin{equation*}
p\Big((m', n',\overline P)|(m,n,l),\omega_k\Big)=[1-(1-\lambda)^{\omega_k}]\hat p_{nn'}\pi^{n'}_{r_{k+1}},
\end{equation*}
where $\hat p_{nn'}$ and $\pi^{n'}_{r_{k+1}}$ are defined in \eqref{eqn:environ_trans} and \eqref{eqn:battery_recur}, respectively.

Similarly, $m'=b_{\rm{max}}$ indicates $r_{k+1}\geqslant b_{\rm{max}}-(m+\omega_k)$ and we have
\begin{align*}
&p\Big(\big(m', n',h(l)\big)|(m,n,l),\omega_k\Big)\\
&=\sum_{r_{k+1}=b_{\rm{max}}-(m+\omega_k)}^{b_{\rm{max}}}(1-\lambda)^{\omega_k}\hat p_{nn'}\pi^{n'}_{r_{k+1}},
\end{align*}
and
\begin{align*}
&p\Big((m', n',\overline P)|(m,n,l), \omega_k\Big)\\
&=\sum_{r_{k+1}=b_{\rm{max}}-(m+\omega_k)}^{b_{\rm{max}}}[1-(1-\lambda)^{\omega_k}]\hat p_{nn'}\pi^{n'}_{r_{k+1}}.
\end{align*}

To formulate \ref{problem:problem-1} into a standard MDP framework, in addition to the state space $\mathbb S$, action set $\mathbb A$ and the one-step state transition probability $\mathbb T: \{p(\phi_2|\phi_1,a)\}$ obtained above, we also need to define the reward functions.

As described in \eqref{eqn:objective-func}, the cost function (objective function) is the trace of average expected state estimate error covariance. Thus we can just define the single stage cost function for time step $k$ as $\mathrm{Tr}\{\mathbb{E}[P_k]\}$, denoted as $v_k(\phi_1,a)$, i.e., as a result of choosing action $\omega_k=a$ when the remote estimator is in state $\Phi_k=\phi_1$ at time step $k$, the remote estimator receive a cost $v_k(\phi,a)$.

Suppose that $v_k(\phi_1, a, \phi_2)$ is the cost given $\Phi_{k+1}=\phi_2=(m',n',l')$, i.e., $P_k=l'$. Thus $v_k(\phi,a)$ can be expressed as the expected value of $v_k(\phi_1, a, \phi_2)$, which depends on the state of the remote estimator at that time step $k$ and at the next time $k+1$:
\begin{align*}
  v_k(\phi_1,a)&=\sum_{\phi_2\in\mathbb S} p( \phi_2| \phi_1, a) v_k(\phi_1, a, \phi_2)\\
  &=(1-\lambda)^a\mathrm{Tr}\{h(l)\}+[1-(1-\lambda)^a]\mathrm{Tr}\{\overline P\},
\end{align*}
where $\Phi_{k}=\phi_1=(m,n,l)$ and $\omega_k=a$.

Without loss of generality, we assume that the costs can be calculated by the sensor prior to selecting a particular action. Define $\Theta$ as the policy for the sensor, which a map from $\mathbb S$ to $\mathbb A$ such that the transmission power is given by $\omega_k=\Theta(\Phi_k)$.

Also denote the expected total cost under a policy $\Theta$ up to time-horizon $T$ when the initial state of the system is $\phi_0$ as
\begin{equation*}
V_T^\Theta(\phi_0)\triangleq \mathbb E_{\phi_0}^\Theta \Big[ \sum_{k=1}^T v_k\big(\phi,a\big)\Big].
\end{equation*}

The performance metric is chosen as the average cost of a policy $\Theta$ given the initial value $\Phi_0=\phi_0$, which is defined by
\begin{equation}\label{eqn:toatl_cost}
J^\Theta(\phi_0)\triangleq \lim_{T\to\infty}\frac{1}{T}V_T^\Theta(\phi_0).
\end{equation}
provided that the limit exists.

\begin{remark}
  Note that if the limit of \eqref{eqn:toatl_cost} does not exist, we can always define
\begin{equation*}
  J^\Theta_-(\phi_0)\triangleq \liminf_{T\to\infty}\frac{1}{T}V_T^\Theta(\phi_0),
\end{equation*}
and
\begin{equation*}
  J^\Theta_+(\phi_0)\triangleq \limsup_{T\to\infty}\frac{1}{T}V_T^\Theta(\phi_0),
\end{equation*}
as the lower and upper bound for $J^\Theta(\phi_0)$ though $J^\Theta_-(\phi_0)$ and $J^\Theta_+(\phi_0)$ may go to infinity.

More detailed stability analysis of $J^\Theta(\phi_0)$ is out of the scope of the current paper and will be left in the future work.
$\hfill \blacksquare$
\end{remark}

Therefore Problem \ref{problem:problem-1} can be stated as finding the optimal policy $\Theta^\star$ to minimize \eqref{eqn:toatl_cost}, i.e.,
\begin{equation*}
  J^\star(\phi_0)=\min_{\Theta} J^\Theta(\phi_0),
\end{equation*}
and
\begin{equation*}
  \Theta^\star=\arg\min_{\Theta} J^\Theta(\phi_0),
\end{equation*}

Based on the theory of MDP, the optimal policy $\Theta^\star$ is stationary and independent of the initial value (\cite{puterman2009markov,bertsekas1995dynamic}). Thus the value of this infinite-time horizon minimization problem is given by $J^\star$ which is the solution of the average-cost optimality (Bellman) equation:
\begin{equation} \label{eqn:bellman}
  J^\star+H(\phi)=\min_{a\in\mathbb A} \Big\{  v_k\big(\phi,a\big) +  \sum_{\phi'\in\mathbb S} p( \phi'| \phi, a) H(\phi')        \Big\},
\end{equation}
where $H$ is the relative value function.

Note that \eqref{eqn:bellman} is not easy to solve (See also \cite{norican13optimal,nayyar2012optimal}). It requires huge computation and cannot be expressed in a closed-form. In addition, as the state set is countably infinite, though the solution can be solved in theory (\cite{bertsekas1995dynamic}), it is quite difficult to implement in practice.  This motivates us to consider a sub-optimal power schedule which can be easily calculated and can be analyzed explicitly.

\section{A Sub-optimal Policy}

In this section, we provide a sub-optimal power schedule policy. In some related literature, the optimal solution is in threshold form (\cite{norican13optimal,nayyar2012optimal,ho2012optimal}), which inspires us to  propose the transmission power schedule  in this  form:
\begin{equation} \label{eqn:sub_opt}
 \omega_k=\left\{\begin{array}{cc}
           \min\{b_k+r_k, R_0\}, & \text{if~} e_k=G, \\
           \min\{b_k+r_k, R_1\}, & \text{if~} e_k=B, \\
         \end{array}
         \right.
\end{equation}
where $R_0(\leqslant b_{\text{max}})$ and $R_1(\leqslant b_{\text{max}})$ are parameters to be designed.

To analyze this strategy, it is convenient to introduce the process $S_k=(b'_k, e_k), k\in\mathbb N,$ for the remote estimator time step $k$, where $b'_k=\min\{b_k+r_k,b_{\text{max}}\}\in\{0,1,2,...,b_{\text{max}}\}$ is the battery level of the sensor after harvesting energy at time step $k$. Based on the description in \eqref{eqn:sub_opt}, it is easy to show that $\{S_k\}$ is a Markov process.

Define the state transition matrix $\Psi=\{\psi_{i,j}\}$, where each element of $\Pi$ is denoted as:
\begin{equation*}
  \psi_{i,j}=\mathbb P[S_{k+1}=(j_1,j_2)|S_k=(i_1,i_2)]
\end{equation*}
where
\begin{align}\label{eqn:mapping}
  i_2&=(i-1)\mod 2,\\
  i_1&=\frac{1}{2}(i-i_2-1),\\
  j_2&=(j-1)\mod 2,\\
  j_1&=\frac{1}{2}(j-j_2-1).
\end{align}

As $i_2,j_2\in\{0,1\}$, it is easy to verify that $i=2i_1+i_2+1$ and $j=2j_1+j_2+1$, thus \eqref{eqn:mapping} is a one-on-one mapping from $\mathbb P[S_{k+1}=(j_1,j_2)|S_k=(i_1,i_2)]$ to $\psi_{i,j}$. Simple analysis leads to the exact form of $\Psi=\{\psi_{i,j}\}$ where:
\begin{align*}
 \psi_{i,j}&=\mathbb P[S_{k+1}=(j_1,j_2)|S_k=(i_1,i_2)]\\
   &=   \left\{\begin{array}{ll}
           \hat p_{i_2j_2}\pi^{j_2}_{j_1}, & \text{if~} i_1<R_{i_2}, \\
           \hat p_{i_2j_2}\pi^{j_2}_{j_1-(i_1-R_{i_2})}, & \text{if~} i_1>R_{i_2} \text{and~} j_1<b_{\text{max}}, \\
           \sum_{m=\overline m}^{b_{\text{max}}}\hat p_{i_2j_2}\pi^{j_2}_{m}, & \text{if~} i_1>R_{i_2} \text{and~} j_1=b_{\text{max}}, \\
         \end{array}
         \right.
\end{align*}
where $\overline m=b_{\text{max}}-(i_1-R_{i_2})$ and $\pi_{j_1}^{j_2}$ is defined in \eqref{eqn:harvested_prob1} and \eqref{eqn:harvested_prob2}.

Here we provide a simple example to illustrate the exact form of $\Psi$. For
example, assume that $b_{\rm{max}}=3$, $R_0=1$, $R_1=2$, then we have $\Psi$ as
in \eqref{eqn:table}, displayed on the following page.

\begin{figure*}[t]
\begin{equation}\label{eqn:table}
\Psi=\left[
  \begin{array}{cccccccccc}
    \hat p_{00}\pi_0^0 &  \hat p_{01}\pi_0^1  &\hat p_{00}\pi_1^0 &  \hat p_{01}\pi_1^1  & \hat p_{00}\pi_2^0 &  \hat p_{01}\pi_2^1  &\hat p_{00}\pi_3^0 &  \hat p_{01}\pi_3^1 \\
    \hat p_{10}\pi_0^0 &  \hat p_{11}\pi_0^1  &\hat p_{10}\pi_1^0 &  \hat p_{11}\pi_1^1  & \hat p_{10}\pi_2^0 &  \hat p_{11}\pi_2^1  &\hat p_{10}\pi_3^0 &  \hat p_{11}\pi_3^1 \\
    \hat p_{00}\pi_0^0 &  \hat p_{01}\pi_0^1  &\hat p_{00}\pi_1^0 &  \hat p_{01}\pi_1^1  & \hat p_{00}\pi_2^0 &  \hat p_{01}\pi_2^1  &\hat p_{00}\pi_3^0 &  \hat p_{01}\pi_3^1 \\
    \hat p_{10}\pi_0^0 &  \hat p_{11}\pi_0^1  &\hat p_{10}\pi_1^0 &  \hat p_{11}\pi_1^1  & \hat p_{10}\pi_2^0 &  \hat p_{11}\pi_2^1  &\hat p_{10}\pi_3^0 &  \hat p_{11}\pi_3^1 \\
    \hat p_{00}\pi_0^0 &  \hat p_{01}\pi_0^1  &\hat p_{00}\pi_1^0 &  \hat p_{01}\pi_1^1  & \hat p_{00}\pi_2^0 &  \hat p_{01}\pi_2^1  &\hat p_{00}\pi_3^0 &  \hat p_{01}\pi_3^1 \\
    0& 0& \hat p_{10}\pi_0^0 &  \hat p_{11}\pi_0^1  &\hat p_{10}\pi_1^0 &  \hat p_{11}\pi_1^1  & \hat p_{10}(\pi_2^0+\pi_3^0)  &  \hat p_{11}(\pi_2^1   +\pi_3^1) \\
    0 & 0 &  \hat p_{00}\pi_0^0 &  \hat p_{01}\pi_0^1  &\hat p_{00}\pi_1^0 &  \hat p_{01}\pi_1^1  & \hat p_{00}(\pi_2^0+\pi_3^0) &  \hat p_{01}(\pi_2^1   +\pi_3^1) \\
    0& 0& 0 & 0 & \hat p_{10}\pi_0^0 &  \hat p_{11}\pi_0^1  & \hat p_{10}(\pi_1^0+\pi_2^0+\pi_3^0)  &  \hat p_{11}(\pi_1^1+\pi_2^1   +\pi_3^1) .
  \end{array}
\right]
\end{equation}
\end{figure*}

Based on the \cite{puterman2009markov,bertsekas1995dynamic}, we can prove that the process described in our work will have a stationary state distribution for each state because this process is a  time-homogeneous Markov chain.

Assume that the stationary state distribution is $q^\star=\{q_0^\star,q_1^\star,q_2^\star,...,q_{2b_{\text{max}}+1}^\star \}$, i.e., in the stationary state,
\begin{equation*}
  \mathbb P[S_k=(i_1,i_2)]=q_{2i_1+i_2}^\star.
\end{equation*}


Based on the power schedule we proposed, $\omega_k$ also has a stationary distribution. Without loss of generality, we assume that $R_0<R_1$. It is easy to derive the stationary distribution for $\omega_k$ :
\begin{equation} \label{eqn:omega_star}
\mathbb P[\omega_k=i]=\left\{
                        \begin{array}{ll}
                          q_{2i}^\star+q_{2i+1}^\star, & \text{if~} 0\leqslant i <R_0, \\
                          \sum_{m=R_0}^{b_{\text{max}}}q_{2m}^\star + q_{2R_0+1}^\star, & \text{if~} i=R_0, \\
                          q_{2i+1}^\star, & \text{if~}  R_0 < i<R_1, \\
                         \sum_{m=R_1}^{b_{\text{max}}}q_{2m+1}^\star , & \text{if~} i=R_1. \\
                          0 , & \text{if~} R_1\leqslant i<b_{\text{max}}. \\
                        \end{array}
                      \right.
\end{equation}

\section{Numerical Simulation}
In this section, we provide a numerical example to illustrate how to implement the sub-optimal solution and evaluate its estimation performance.

Consider a scalar system  with parameters $A=0.9, C=0.7, R=Q=0.8, \lambda=0.7$.  Note that even for scalar systems, the states set for the process is still infinite, which renders finding the optimal solution intractable. Thus we will compare the suboptimal one with other policies. Assume that $b_{\rm{max}}=3$, $R_0=1$, and $R_1=2$. The environment condition transition probabilities are set as $\hat p_{00}=0.7, \hat p_{01}=0.3, \hat p_{10}=0.2, \hat p_{11}=0.8$ and the distribution of harvested energy is defined as
\begin{equation*}
\pi_i^0=\left\{
            \begin{array}{cc}
            0.1, & i=0,\\
            0.2, & i=1,\\
            0.3, & i=2,\\
            0.4, & i=3,\\
            \end{array}
        \right.
\end{equation*}
and
\begin{equation*}
\pi_i^1=\left\{
            \begin{array}{cc}
            0.4, & i=0,\\
            0.3, & i=1,\\
            0.2, & i=2,\\
            0.1, & i=3.\\
            \end{array}
        \right.
\end{equation*}

Based on \eqref{eqn:table}, we can easily calculate $\Psi$:
\begin{equation*}
\Psi=\left[
  \begin{array}{cccccccccc}
    0.07 & 0.12 & 0.14 & 0.09 & 0.21 & 0.06 & 0.28 & 0.03\\
    0.02 & 0.32 & 0.04 & 0.24 & 0.06 & 0.16 & 0.08 & 0.08\\
    0.07 & 0.12 & 0.14 & 0.09 & 0.21 & 0.06 & 0.28 & 0.03\\
    0.02 & 0.32 & 0.04 & 0.24 & 0.06 & 0.16 & 0.08 & 0.08\\
    0.07 & 0.12 & 0.14 & 0.09 & 0.21 & 0.06 & 0.28 & 0.03\\
    0     &  0    & 0.02  & 0.32  & 0.04 & 0.24 & 0.14 & 0.24\\
    0     & 0     & 0.07 & 0.12 & 0.14 & 0.09 & 0.49 & 0.09\\
    0     & 0     & 0      &    0  & 0.02 & 0.32 & 0.18 & 0.48
  \end{array}
\right]
\end{equation*}
and obtain the stationary distribution $q^\star $.

Define
\begin{equation*}
  J_k(\theta)=\frac{1}{k}\sum_{i=1}^{k}\mathrm{Tr} \left(\mathbb{E}[P_i]\right),
\end{equation*}
as the empirical approximation (via 100000 Monte Carlo simulations) of $J(\theta)$ (See \eqref{eqn:objective-func}) at every time instant $k$.

\begin{figure}[t]

  \centering

  \includegraphics[width=9cm]{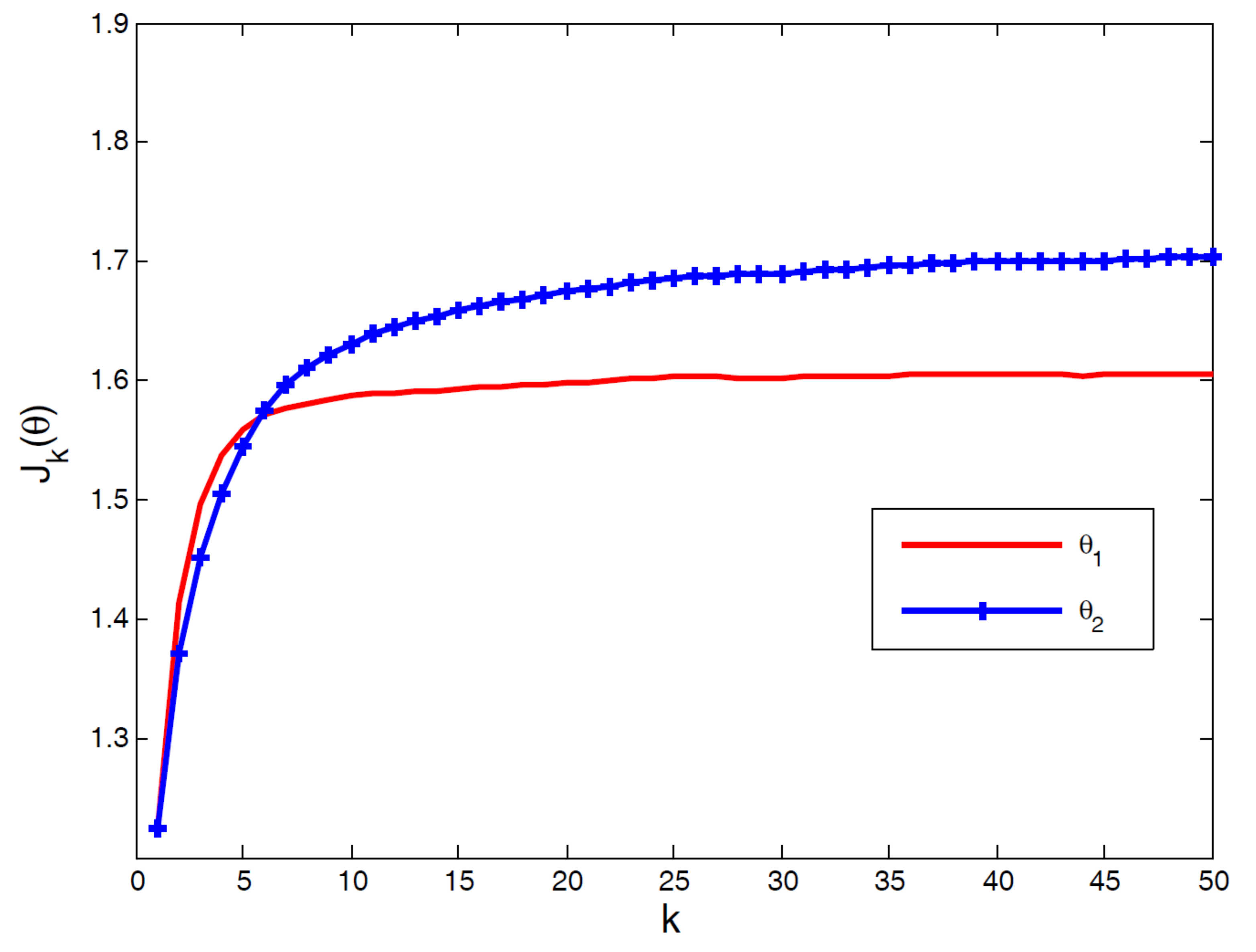}

  \caption{Estimation performance comparison of $\theta_1$ and $\theta_2$} \label{fig:compare1}

\end{figure}

As a comparison, we propose another common transmission power schedule, i.e, the ``greedy'' method:
\begin{equation*}
  \omega_k=r_k,
\end{equation*}
which refers to using all the harvested energy $r_k$ to send the data packet at each time step. Denote our proposed sub-optimal schedule as $\theta_1$, and the ``greedy'' method as $\theta_2$. Though both methods are easy to implement, the simulation shows that our proposed sub-optimal method $\theta_1$ obtains a better estimation performance. (See Fig. \ref{fig:compare1} ). Note that the ``greedy'' method have a better performance only in the first several time steps, which is because the ``greedy'' method used all the harvested energy instead of reserve some for the future. The more cautious energy management policy \eqref{eqn:sub_opt} makes better use of the battery capability.

\section{Conclusion}
We have studied remote estimation with  a wireless sensor in this paper. Instead of using a conventional battery-powered sensor, a sensor equipped with an energy harvester which can obtain energy from the external environment was utilized. We formulated this problem into an infinite time-horizon Markov decision process and provide the optimal sensor transmission power control strategy. In addition, a sub-optimal policy which is easier to implement and requires less computations is also presented. Numerical simulations illustrate that performance gains can be obtained when compared to a greedy method.

\bibliographystyle{ieeetr}
\bibliography{reference1}

\begin{thebibliography}{10}

\bibitem{yick2008wireless}
J.~Yick, B.~Mukherjee, and D.~Ghosal, ``Wireless sensor network survey,'' {\em
  Computer networks}, vol.~52, no.~12, pp.~2292--2330, 2008.

\bibitem{ho2012optimal}
C.~K. Ho and R.~Zhang, ``Optimal energy allocation for wireless communications
  with energy harvesting constraints,'' {\em IEEE Transactions on Signal
  Processing}, vol.~60, no.~9, pp.~4808--4818, 2012.

\bibitem{aziz2013survey}
A.~A. Aziz, Y.~A. Sekercioglu, P.~Fitzpatrick, and M.~Ivanovich, ``A survey on
  distributed topology control techniques for extending the lifetime of battery
  powered wireless sensor networks,'' {\em IEEE Communications Surveys and
  Tutorials}, vol.~15, no.~1, pp.~121--144, 2013.

\bibitem{pantazis2007survey}
N.~A. Pantazis and D.~D. Vergados, ``A survey on power control issues in
  wireless sensor networks,'' {\em IEEE Communications Surveys and Tutorials},
  vol.~9, no.~4, pp.~86--107, 2007.

\bibitem{quevedo2010energy}
D.~E. Quevedo, A.~Ahl{\'e}n, and J.~{\O}stergaard, ``Energy efficient state
  estimation with wireless sensors through the use of predictive power control
  and coding,'' {\em IEEE Transactions on Signal Processing}, vol.~58, no.~9,
  pp.~4811--4823, 2010.

\bibitem{sudevalayam2011energy}
S.~Sudevalayam and P.~Kulkarni, ``Energy harvesting sensor nodes: Survey and
  implications,'' {\em IEEE Communications Surveys and Tutorials}, vol.~13,
  no.~3, pp.~443--461, 2011.

\bibitem{nayyar2012optimal}
A.~Nayyar, T.~Basar, D.~Teneketzis, and V.~V. Veeravalli, ``Optimal strategies
  for communication and remote estimation with an energy harvesting sensor,''
  {\em IEEE Transactions on Automatic Control}, vol.~58, no.~9, pp.~2246--2260,
  2013.

\bibitem{li2013optimal}
Y.~Li, D.~E. Quevedo, V.~Lau, and L.~Shi, ``Optimal periodic transmission power
  schedules for remote estimation of {ARMA} processes,'' {\em IEEE Transactions
  on Signal Processing}, vol.~61, no.~24, pp.~6164--6174, 2013.

\bibitem{norican13optimal}
M.~Nourian, A.~Leong, and S.~Dey, ``Optimal energy allocation for {Kalman}
  filtering over packet dropping links with energy harvesting constraints,'' in
  {\em 4th IFAC Workshop on Distributed Estimation and Control in Networked
  Systems, Koblenz, Germany}, 2013.

\bibitem{hovareshti2007sensor}
P.~Hovareshti, V.~Gupta, and J.~S. Baras, ``Sensor scheduling using smart
  sensors,'' in {\em {Proceedings of the 46th IEEE Conference on Decision and
  Control}}, pp.~494--499, 2007.

\bibitem{shi2012scheduling}
L.~Shi and H.~Zhang, ``{Scheduling two Gauss-Markov systems: an optimal
  solution for remote state estimation under bandwidth constraint},'' {\em IEEE
  Transactions on Signal Processing}, vol.~60, no.~4, pp.~2038--2042, 2012.

\bibitem{wu2013can}
J.~Wu, Y.~Yuan, H.~Zhang, and L.~Shi, ``How can online schedules improve
  communication and estimation tradeoff?,'' {\em IEEE Transactions on Signal
  Processing}, vol.~61, no.~7, pp.~1625--1631, 2013.

\bibitem{li2013jamming}
Y.~Li, L.~Shi, P.~Cheng, J.~Chen, and D.~E. Quevedo, ``Jamming attack on
  cyber-physical systems: A game-theoretic approach,'' in {\em IEEE
  International Conference on CYBER Technology in Automatation, Control, and
  Intelligent Systems, Nanjing, China}, 2013.

\bibitem{anderson1981detectability}
B.~D.~O. Anderson and J.~B. Moore, ``Detectability and stabilizability of
  time-varying discrete-time linear systems,'' {\em SIAM Journal on Control and
  Optimization}, vol.~19, no.~1, pp.~20--32, 1981.

\bibitem{shi2011time}
L.~Shi, K.~H. Johansson, and L.~Qiu, ``Time and event-based sensor scheduling
  for networks with limited communication resources,'' in {\em World Congress
  of the International Federation of Automatic Control (IFAC)}, vol.~18,
  pp.~13263--13268, 2011.

\bibitem{li2013online}
Y.~Li, D.~E. Quevedo, V.~Lau, and L.~Shi, ``Online sensor transmission power
  schedule for remote state estimation,'' in {\em IEEE 52nd Aunnal Conference
  on Decision and Control (CDC), Florence, Italy}, 2013.

\bibitem{shi2010kalman}
L.~Shi, M.~Epstein, and R.~M. Murray, ``Kalman filtering over a packet-dropping
  network: A probabilistic perspective,'' {\em IEEE Transactions on Automatic
  Control}, vol.~55, no.~3, pp.~594--604, 2010.

\bibitem{puterman2009markov}
M.~L. Puterman, {\em Markov Decision Processes: Discrete Stochastic Dynamic
  Programming}, vol.~414.
\newblock Wiley. com, 2009.

\bibitem{bertsekas1995dynamic}
D.~P. Bertsekas, {\em Dynamic Programming and Optimal Control}.
\newblock Athena Scientific Belmont, 1995.

\end{thebibliography}

\end{document}